\documentclass[12pt]{article}

\title{
Nonparametric estimation in a  semimartingale regression model.\\
Part 2. Robust asymptotic efficiency.
\thanks{The paper is  supported by the RFFI-Grant 09-01-00172-a.}
}

\author{Konev, V.
\thanks{
Department of Applied Mathematics and Cybernetics,
 Tomsk State University,
Lenin str. 36,
 634050 Tomsk, Russia,
 e-mail: vvkonev@mail.tsu.ru }
 \and
S. Pergamenshchikov\thanks{
 Laboratoire de Math\'ematiques Raphael Salem,
 Avenue de l'Universit\'e, BP. 12,
  Universit\'e de Rouen,
   F76801, Saint Etienne du Rouvray, Cedex France,
   Department of Mathematics and Mechanics,Tomsk State University,
Lenin str. 36, 634041 Tomsk, Russia,
 e-mail: Serge.Pergamenchtchikov@univ-rouen.fr}
}

\usepackage{amssymb}
\usepackage{amsfonts}
\usepackage{amsmath}
\usepackage{amsthm}
\usepackage{enumerate}
\usepackage{multicol}

\newtheorem{theorem}{Theorem}[section]
\newtheorem{proposition}[theorem]{Proposition}

\newtheorem{remark}{Remark}[section]
\newtheorem{corollary}[theorem]{Corollary}

\newcommand\cA{{\cal A}}
\newcommand\cC{{\cal C}}

\newcommand\cH{{\cal H}}

\newcommand\cL{{\cal L}}
\newcommand\cB{{\cal B}}

\newcommand\cN{{\cal N}}

\newcommand\cP{{\cal P}}
\newcommand\cQ{{\cal Q}}
\newcommand\cX{{\cal X}}
\newcommand\cD{{\cal D}}

\newcommand\cR{{\cal R}}

\newcommand\ve{\varepsilon}
\newcommand\ov{\overline}

\def\bbr{{\mathbb R}}

\def\text#1{\hbox{#1}}
\def\proof{{\noindent \bf Proof. }}
\def\endproof{\mbox{\ $\qed$}}

\def\E{{\bf E}}
\def\B{{\bf B}}

\def\C{{\bf C}}

\def\L{{\bf L}}

\newcommand{\wh}{\widehat}
\newcommand{\wt}{\widetilde}

\def\Chi{{\bf 1}}
\def\d{\mathrm{d}}
\def\build #1_#2{\mathrel{\mathop{\kern 0pt #1}\limits_{#2}}}
\newcommand{\zs}[1]{{\mathchoice{#1}{#1}{\lower.25ex\hbox{$\scriptstyle#1$}}
{\lower0.25ex\hbox{$\scriptscriptstyle#1$}}}}
\numberwithin{equation}{section}

\begin{document}

\maketitle

\begin{abstract}
In this paper we prove the asymptotic efficiency of the model
selection procedure proposed by the authors in
\cite{KoPe3}. To
this end we introduce the robust risk as the least
upper bound of the quadratical risk over a broad class of
observation distributions.  Asymptotic upper and lower bounds
 for the robust risk have been derived. The asymptotic efficiency
 of the procedure is proved. The Pinsker constant is found.
\end{abstract}
\vspace*{5mm}

\vspace*{5mm}

\bibliographystyle{plain}
\renewcommand{\columnseprule}{.1pt}

\newpage

\section{Introduction}\label{sec:In}

In this paper we will investigate the asymptotic efficiency of the
model selection procedure proposed in \cite{KoPe3} for estimating
a 1-periodic function $S\,:\,\bbr\to\bbr$, $S\in\cL_\zs{2}[0,1]$,
in a continuous time regression model
 \begin{equation}\label{sec:In.1}
  \d y_\zs{t}=S(t)\d t+\d \xi_t\,,
  \quad 0\le t\le n\,,
 \end{equation}
with a semimartingale noise $\xi=(\xi_\zs{t})_\zs{0\le t\le n}$.
The quality of an estimate $\wt{S}$ (any real-valued function measurable with
respect to $\sigma\{y_\zs{t}\,,\,0\le t\le n\}$) for $S$ is given by the mean
integrated squared error, i.e.
\begin{equation}\label{sec:In.5}
\cR_\zs{Q}(\wt{S},S) = \E_\zs{Q,S}\,||\wt{S}-S||^{2}\,,
\end{equation}
where  $\E_\zs{Q,S}$ is the expectation with respect to the
noise distribution $Q$ given a function $S$;
$$
||S||^{2}=\int^{1}_\zs{0}\,S^{2}(x)\d x\,.
$$
The semimartingale noise $(\xi_\zs{t})_\zs{0\le t\le n}$ is assumed to take values in the
Skorohod space $\cD[0,n]$ and has the distribution $Q$ on $\cD[0,n]$
 such that for any function $f$ from $\cL_\zs{2}[0,n]$ the stochastic integral
\begin{equation}\label{sec:In.5-0}
I_\zs{n}(f)=\int^n_\zs{0}f_\zs{s}\d \xi_\zs{s}
\end{equation}
is well defined with
\begin{equation}\label{sec:In.5-1}
\E_\zs{Q} I_\zs{n}(f)=0
\quad\mbox{and}\quad
\E_\zs{Q} I^2_\zs{n}(f)\le \sigma^{*}\,\int^n_\zs{0}\,f^2_\zs{s}\,\d s
\end{equation}
where $\sigma^*$ is some positive constant which may, in general, depend on $n$, i.e.
$\sigma^{*}=\sigma^{*}_\zs{n}$, such that
\begin{equation}\label{sec:In.12-1}
0<\liminf_\zs{n\to\infty}\sigma^{*}_\zs{n}\le
\limsup_\zs{n\to\infty}\sigma^{*}_\zs{n}<\infty
\,.
\end{equation}

Now we  define
a robust risk function which is required to measure the quality of
an estimate $\wt{S}$  provided that  a true
distribution of the noise $(\xi_\zs{t})_\zs{0\le t\le n}$ is known
 to  belong to some family of distributions $\cQ^{*}_\zs{n}$
which will be specified below. Just as in  \cite{GaPe3}  we define
the robust risk as
\begin{equation}\label{sec:In.12}
\cR^{*}_\zs{n}(\wt{S}_\zs{n},S)=\sup_\zs{Q\in\cQ^{*}_\zs{n}}\,
\cR_\zs{Q}(\wt{S}_\zs{n},S)\,.
\end{equation}

The goal of this paper is to prove that the model selection procedure for estimating $S$
in the model \eqref{sec:In.1}
constructed in \cite{KoPe3}  is asymptotically efficient with respect to this risk.
When studying the asymptotic efficiency of this procedure, described in detail
in Section~\ref{sec:Ep},
 we suppose that the unknown function $S$ in the model
 \eqref{sec:In.1} belongs to the Sobolev ball
\begin{equation}\label{sec:In.13}
W^{k}_\zs{r}=\{f\in \,\cC^{k}_\zs{per}[0,1]
\,,\,\sum_\zs{j=0}^k\,||f^{(j)}||^2\le r\}\,,
 \end{equation}
where $r>0\,,\ k\ge 1$ are some  parameters,
$\cC^{k}_\zs{per}[0,1]$ is a set of
 $k$ times continuously differentiable functions
$f\,:\,[0,1]\to\bbr$ such that $f^{(i)}(0)=f^{(i)}(1)$ for all
$0\le i \le k$. The functional class $W^{k}_\zs{r}$ can be written as
the ellipsoid in $l_\zs{2}$, i.e.
 \begin{equation}\label{sec:In.14}
W^{k}_\zs{r}=\{f\in\,\cC^{k}_\zs{per}[0,1]\,:\,
\sum_\zs{j=1}^{\infty}\,a_\zs{j}\,\theta^2_\zs{j}\,\le r\}
 \end{equation}
where
$$
a_\zs{j}=\sum^k_{i=0}\left(2\pi [j/2]\right)^{2i}\,.
$$

In \cite{KoPe3} we established a sharp non-asymptotic oracle inequality
for mean integrated squared error \eqref{sec:In.5}. The proof of the asymptotic
efficiency of the model selection procedure below largely bases
on the counterpart of this inequality for the robust risk \eqref{sec:In.12}
given in Theorem~\ref{Th.sec:Or.1}.

 It will be observed that the notion "nonparametric robust risk" was initially introduced
 in \cite{GaPe0} for estimating a regression curve at a fixed point. The greatest
 lower bound for such risks have been derived and a point estimate is found for which
this bound is attained. The latter means that the point estimate turns out to
be robust efficient. In \cite{Br} this approach was applied for pointwise estimation
in a heteroscedastic regression model.

The optimal convergence rate of the robust quadratic risks has been obtained in
\cite{KoPe2} for the non-parametric estimation problem in a continuous time
regression model with a coloured noise having unknown correlation properties
under full and partial observations. The asymptotic efficiency with respect to
the robust quadratic risks, has been studied in \cite{GaPe3}, \cite{GaPe4}
for the problem of non-parametric estimation in heteroscedastic regression
models. In this paper we apply this approach for the model \eqref{sec:In.1}.

The rest of the paper is organized as follows. In Section~\ref{sec:Ep}
we construct the model selection procedure and
formulate (Theorem 2.1) the oracle inequality for
the robust risk. Section~\ref{sec:Mr} gives the main results.
In Section~\ref{sec:Ex} we consider an example of the model \eqref{sec:In.1} with
 the Levy type  martingale noise. In Section \ref{sec:Up} and \ref{sec:Lo}
 we obtain the upper and lower bounds for the robust risk.
In Section~\ref{sec:A} some technical results are established.

\section{Oracle inequality for the robust risk }\label{sec:Ep}
The model selection procedure is constructed on the basis of a
weighted least squares estimate having the form
\begin{equation}\label{sec:In.1-1}
\wh{S}_\zs{\gamma}=\sum^{\infty}_\zs{j=1}\gamma(j)\wh{\theta}_\zs{j,n}\phi_\zs{j}
\quad\mbox{with}\quad \wh{\theta}_\zs{j,n}=
\frac{1}{n}\int^n_\zs{0}\,\phi_j(t)\,\d y_\zs{t}\,,
\end{equation}
where $(\phi_\zs{j})_\zs{j\ge 1}$ is the standard trigonometric
basis in $\cL_\zs{2}[0,1]$ defined as
\begin{equation}\label{sec:In.2}
\phi_1=1\,,\quad \phi_\zs{j}(x)=\sqrt{2}\,Tr_\zs{j}(2\pi
[j/2]x)\,,\ j\ge 2\,,
\end{equation}
where the function $Tr_\zs{j}(x)=\cos(x)$ for even $j$ and
$Tr_\zs{j}(x)=\sin(x)$ for odd $j$; $[x]$ denotes the integer part
of $x$. The sample functionals $\wh{\theta}_\zs{j,n}$ are estimates
of the corresponding Fourier coefficients
\begin{equation}\label{sec:In.3}
\theta_\zs{j}=(S,\phi_j)= \int^1_\zs{0}\,S(t)\,\phi_\zs{j}(t)\,\d
t\,.
\end{equation}
Further we introduce the cost function as
$$
J_\zs{n}(\gamma)\,=\,\sum^\infty_\zs{j=1}\,\gamma^2(j)\wh{\theta}^2_\zs{j,n}\,-
2\,\sum^\infty_\zs{j=1}\,\gamma(j)\,\wt{\theta}_\zs{j,n}\,
+\,\rho\,\wh{P}_\zs{n}(\gamma)\,.
$$
Here
$$
\wt{\theta}_\zs{j,n}=
\wh{\theta}^2_\zs{j,n}-\frac{\wh{\sigma}_\zs{n}}{n}
\quad\mbox{with}\quad
\wh{\sigma}_\zs{n}=\sum^n_\zs{j=l}\,\wh{\theta}^2_\zs{j,n}\,,
\quad l=[\sqrt{n}]+1\,;
$$
$\wh{P}_\zs{n}(\gamma)$ is the penalty term defined as
$$
\wh{P}_\zs{n}(\gamma)=\frac{\wh{\sigma}_\zs{n}\,|\gamma|^2}{n} \,.
$$
As to the parameter $\rho$, we assume that this parameter is a function of
$n$, i.e. $\rho=\rho_\zs{n}$ such that
 $0<\rho<1/3$ and
$$
\lim_\zs{n\to\infty}\,n^{\delta}\,\rho_\zs{n}=0
\quad\mbox{for all}\quad
\delta>0\,.
$$
We define the model selection procedure as
\begin{equation}\label{sec:In.4}
\wh{S}_\zs{*}=\wh{S}_\zs{\wh{\gamma}}
\end{equation}
where $\wh{\gamma}$ is the minimizer of the cost function $J_\zs{n}(\gamma)$
in some given class $\Gamma$ of weight sequences
$\gamma=(\gamma(j))_\zs{j\ge 1}\in [0,1]^{\infty}$, i.e.
\begin{equation}\label{sec:In.4-1}
\wh{\gamma}=\mbox{argmin}_\zs{\gamma\in\Gamma}\,J_n(\gamma)\,.
\end{equation}

Now we specify the family of distributions $\cQ^{*}_\zs{n}$ in the robust risk
\eqref{sec:In.12}. Let $\cP_\zs{n}$ denote the class of all  distributions $Q$
of the semimartingale $(\xi_\zs{t})$ satisfying the condition
\eqref{sec:In.5-1}.
It is obvious that the distribution $Q_\zs{0}$ of the process
$\xi_\zs{t}=\sqrt{\sigma^{*}}w_\zs{t}$, where $(w_\zs{t})$
is a standard Brownian motion, enters  the class $\cP_\zs{n}$, i.e.
$Q\in\cP_\zs{n}$. In addition, we need to impose some technical conditions
on the distribution $Q$ of the process $(\xi_\zs{t})_\zs{0\le t\le n}$.
Let denote
\begin{equation}\label{sec:In.6}
\sigma(Q)=\ov{\lim}_\zs{n\to\infty}\max_\zs{1\le j\le n}\,\E_\zs{Q}\,\xi^{2}_\zs{j,n}\,,
\end{equation}
where
$$
\xi_\zs{j,n}=\frac{1}{\sqrt{n}}\,I_\zs{n}(\phi_\zs{j})\,,
$$
($I_\zs{n}(\phi_\zs{j})$
is given in \eqref{sec:In.5-0}) and introduce
 two $\cP_\zs{n}\to\bbr_\zs{+}$ functionals
$$
\L_\zs{1,n}(Q)=
 \sup_\zs{x\in \cH\,,\,\#(x)\le n}\,
 \left|
\sum^\infty_\zs{j=1}\,x_\zs{j}\,
\left(
\E_\zs{Q}\,\xi^{2}_\zs{j,n}
-
\sigma(Q)
\right)
\right|
$$
and
$$
\L_\zs{2,n}(Q)=
 \sup_\zs{|x|\le 1\,,\,\#(x)\le n}\,
\E_\zs{Q}\, \left(
\sum^\infty_\zs{j=1}\,x_\zs{j}\,
\wt{\xi}_\zs{j,n}\,
\right)^2
$$
where $\cH=[-1,1]^{\infty}$, $|x|^2=\sum^\infty_\zs{j=1}x^2_\zs{j}$,
$\#(x)=\sum^\infty_\zs{j=1}\,\Chi_\zs{\{|x_\zs{j}|>0\}}$ and
$$
\wt{\xi}_\zs{j,n}=\xi^{2}_\zs{j,n}-\E_\zs{Q}\xi^{2}_\zs{j,n}\,.
$$
Now we consider the family of all distributions $Q$ from $\cP_\zs{n}$
with the growth restriction on $\L_\zs{1,n}(Q)+\L_\zs{2,n}(Q)$, i.e.
$$
\cP^{*}_\zs{n}=\left\{Q\in\cP_\zs{n}\,:\, \L_\zs{1,n}(Q)+\L_\zs{2,n}(Q)\le
\,l_\zs{n} \right\}\,,
$$
where $l_\zs{n}$ is a slowly increasing positive function, i.e.
$l_\zs{n}\to+\infty$ as $n\to+\infty$ and for any $\delta>0$
$$
\lim_\zs{n\to\infty}\frac{l_\zs{n}}{n^{\delta}}=0\,.
$$
It will be observed that any distribution $Q$ from $\cP^{*}_\zs{n}$
satisfies conditions $\C_\zs{1})$ and $\C_\zs{2})$ on the noise distribution
from  \cite{KoPe3} with $c^{*}_\zs{1,n}\le l_\zs{n}$ and
$c^{*}_\zs{2,n}\le l_\zs{n}$. We remind that these conditions are

\noindent $\C_\zs{1})$
$$
c^{*}_\zs{1,n}=\L_\zs{1,n}(Q) < \infty\,;
$$
\noindent $\C_\zs{2})$
$$
c^{*}_\zs{2,n}=\L_\zs{2,n}(Q) < \infty\,.
$$

\noindent

In the sequel we assume that the distribution of the noise $(\xi_\zs{t})$
in \eqref{sec:In.1} is known up to its belonging to some
distribution family satisfying the following condition.\\[2mm]

\noindent $\C^{*})$ {\sl
Let $\cQ^{*}_\zs{n}$ be a family of the distributions $Q$ from
$\cP^{*}_\zs{n}$ such that $Q_\zs{0}\in\cQ^{*}_\zs{n}$.
}

\vspace{4mm}

\noindent An important example for such family is given in Section~\ref{sec:Ex}.

Now we specify the set $\Gamma$ in the model
selection procedure \eqref{sec:In.4} and state the oracle inequality
for the robust risk \eqref{sec:In.12}
which is a counterpart of that
obtained in  \cite{KoPe3} for the mean integrated squared error \eqref{sec:In.5}.
 Consider the numerical grid
\begin{equation}\label{sec:Or.1}
\cA_\zs{n}=\{1,\ldots,k^*\}\times\{t_1,\ldots,t_m\}\,,
\end{equation}
 where
$t_i=i\ve$ and $m=[1/\ve^2]$; parameters
$k^*\ge 1$ and $0<\ve\le 1$ are functions of $n$, i.e. $
k^*=k^*(n)$ and $\ve=\ve(n)$, such that
for any $\delta>0$
\begin{equation}\label{sec:Or.1-1}
\left\{
\begin{array}{ll}
&\lim_\zs{n\to\infty}\,k^*(n)=+\infty\,,
\quad
\lim_\zs{n\to\infty}\,\dfrac{k^*(n)}{\ln n}=0\,,\\[6mm]
&\lim_\zs{n\to\infty}\ve(n)=0
\quad\mbox{and}\quad
\lim_\zs{n\to\infty}\,n^{\delta}\ve(n)\,=+\infty\,.
\end{array}
\right.
\end{equation}
 For example, one can take
$$
\ve(n)=\frac{1}{\ln (n+1)}
\quad\mbox{and}\quad
k^*(n)=\sqrt{\ln (n+1)}
$$
for  $n\ge 1$.

Define the set $\Gamma$ as
\begin{equation}\label{sec:Or.2}
\Gamma\,=\,\{\gamma_\zs{\alpha}\,,\,\alpha\in\cA_\zs{n}\}\,,
\end{equation}
where $\gamma_\zs{\alpha}$ is the weight sequence corresponding to an element
$\alpha=(\beta,t)\in\cA_\zs{n}$, given
by the formula
\begin{equation}\label{sec:Or.2-1}
\gamma_\zs{\alpha}(j)=\Chi_\zs{\{1\le j\le j_\zs{0}\}}+
\left(1-(j/\omega_\alpha)^\beta\right)\, \Chi_\zs{\{ j_\zs{0}<j\le
\omega_\alpha\}}
\end{equation}
where $j_\zs{0}=j_\zs{0}(\alpha)=\left[\omega_\zs{\alpha}/(1+\ln n)\right]$,
$\omega_\zs{\alpha}=(\tau_\zs{\beta}\,t\,n)^{1/(2\beta+1)}$ and
$$
\tau_\zs{\beta}=\frac{(\beta+1)(2\beta+1)}{\pi^{2\beta}\beta}\,.
$$
Along the lines of the proof of Theorem 2.1 in \cite{KoPe3} one can
establish the following result.

\begin{theorem}\label{Th.sec:Or.1}
Assume that the unknown function $S$ is continuously differentiable
and the distribution family $\cQ^{*}_\zs{n}$ in the robust
 risk \eqref{sec:In.12} satisfies the condition $\C^{*})$.
Then
   the estimator
\eqref{sec:In.4},
for any $n\ge 1$,
 satisfies the oracle inequality
\begin{align}\label{sec:Or.3}
\cR^{*}_\zs{n}(\wh{S}_\zs{*},S)\,\le\, \frac{1+3\rho-2\rho^2}{1-3\rho}
\min_\zs{\gamma\in\Gamma} \cR^{*}_\zs{n}(\wh{S}_\zs{\gamma},S)
+\frac{1}{n}\,\cD_\zs{n}(\rho)\,,
\end{align}
where the term $\cD_\zs{n}(\rho)$ is defined in \cite{KoPe3} such that
\begin{equation} \label{sec:Or.4}
\lim_\zs{n\to\infty}\,\frac{\cD_\zs{n}(\rho)}{n^{\delta}}=0
\end{equation}
for each $\delta>0$.
\end{theorem}

\begin{remark}\label{Re.sec:Or.1}
The inequality \eqref{sec:Or.3} will be used to derive
the upper bound for the robust risk  \eqref{sec:In.12}. It will be noted
that the second summand in \eqref{sec:Or.3} when multiplied by the optimal rate
$n^{2k/(2k+1)}$ tends to zero as $n\to\infty$ for each $k\ge 1$. Therefore, taking into
account that $\rho\to 0$ as $n\to\infty$, the principal term in the upper bound
is given by the minimal risk over the family of estimates
$(\wh{S}_\zs{\gamma})_\zs{\gamma\in\Gamma}$.
As is shown in \cite{GaPe2},
 the efficient estimate enters this family. However one can not
use this estimate because it depends on the unknown parameters $k\ge 1$ and $r>0$ of the Sobolev
ball. It is this fact that shows an adaptive role of the oracle inequality
\eqref{sec:Or.3} which gives the asymptotic upper bound in the case
when this information is not available.
\end{remark}

\section{Main results}\label{sec:Mr}

In this Section we will show, proceeding from \eqref{sec:Or.3}, that
the Pinsker constant for the robust risk \eqref{sec:In.12}
is given by the equation
\begin{equation}\label{sec:Mr.1}
R^{*}_\zs{k,n}=\left((2k+1)r\right)^{1/(2k+1)}\,
\left(
\frac{\sigma^{*}_\zs{n}k}{(k+1)\pi} \right)^{2k/(2k+1)}\,.
\end{equation}
It is well known that
 the optimal (minimax) rate
for the Sobolev ball $W^{k}_\zs{r}$
 is $n^{2k/(2k+1)}$ (see, for example, \cite{Pi}, \cite{Nu}).
We will see that asymptotically the robust risk
of the model selection \eqref{sec:In.4} normalized by this rate is bounded
from above by $R^{*}_\zs{k,n}$. Moreover, this bound can not be diminished if one
considers the class of all admissible estimates for $S$.

\begin{theorem}\label{Th.sec:Mr.1}
Assume that, in model \eqref{sec:In.1}, the distribution
of $(\xi_\zs{t})$ satisfies the
 condition $\C^{*})$.
 Then  the
robust risk \eqref{sec:In.12}
of the  model selection estimator $\wh{S}_\zs{*}$ defined in \eqref{sec:In.4}, \eqref{sec:Or.2},
 has
 the following asymptotic upper bound
 \begin{equation}\label{sec:Mr.2}
\limsup_\zs{n\to\infty}\,n^{2k/(2k+1)}\,
\frac{1}{R^{*}_\zs{k,n}}
 \sup_\zs{S\in W^k_r}\,
\cR^{*}_\zs{n}(\wh{S}_\zs{*},S) \le  1\,.
 \end{equation}
\end{theorem}
Now we obtain a lower bound for the robust risk
\eqref{sec:In.12}.
Let $\Pi_\zs{n}$ be the set of all estimators $\wt{S}_\zs{n}$
measurable with respect to the sigma-algebra

$\sigma\{y_\zs{t}\,,\,0\le t\le n\}$
 generated by the process \eqref{sec:In.1}.
\begin{theorem}\label{Th.sec:Mr.2} Under the conditions of Theorem~\ref{Th.sec:Mr.1}
 \begin{equation}\label{sec:Mr.3}
\liminf_{n\to\infty}\,n^{2k/(2k+1)}\,\frac{1}{R^{*}_\zs{k,n}}\, \inf_{\wt{S}_\zs{n}\in\Pi_\zs{n}}\,\,
\sup_\zs{S\in W^{k}_\zs{r}} \,\cR^{*}_\zs{n}(\wt{S}_\zs{n},S) \ge
1 \,.
 \end{equation}
\end{theorem}
\noindent Theorem~\ref{Th.sec:Mr.1} and Theorem~\ref{Th.sec:Mr.2} imply the following result
\begin{corollary}\label{Co.sec:Mr.1}
Under the conditions of Theorem~\ref{Th.sec:Mr.1}
\begin{equation}\label{sec:Mr.4}
\lim_{n\to\infty}\,n^{2k/(2k+1)}\,\frac{1}{R^{*}_\zs{k,n}}\, \inf_{\wt{S}_\zs{n}\in\Pi_\zs{n}}\,\,
\sup_\zs{S\in W^{k}_\zs{r}} \,\cR^{*}_\zs{n}(\wt{S}_\zs{n},S)
=1 \,.
 \end{equation}
\end{corollary}

\begin{remark}\label{Re.sec:Mr.1}
The equation \eqref{sec:Mr.4} means that  the sequence $R^{*}_\zs{k,n}$ defined by \eqref{sec:Mr.1}
is the Pinsker constant (see, for example, \cite{Pi}, \cite{Nu}) for the model \eqref{sec:In.1}.
\end{remark}

\section{Example}\label{sec:Ex}

Let the process $(\xi_\zs{t})$ be defined as
\begin{equation}\label{sec:Ex.1}
\xi_\zs{t}=\varrho_\zs{1}w_\zs{t}+\varrho_\zs{2}z_\zs{t}\,,
\end{equation}
where $(w_\zs{t})_\zs{t\ge 0}$ is a standard Brownian motion, $(z_\zs{t})_\zs{t\ge 0}$
is a compound Poisson process defined as
$$
z_\zs{t}=\sum^{N_\zs{t}}_\zs{j=1}\,Y_\zs{j}\,,
$$
where $(N_\zs{t})_\zs{t\ge 0}$ is a standard homogeneous Poisson process with
unknown  intensity $\lambda>0$ and
$(Y_\zs{j})_\zs{j\ge 1}$ is an i.i.d. sequence of random variables with
$$
\E\, Y_\zs{j}=0\,,
\quad
\E\, Y^2_\zs{j}=1
\quad\mbox{and}\quad
\E\, Y^{4}_\zs{j}\,<\infty\,.
$$
Substituting
\eqref{sec:Ex.1} in
\eqref{sec:In.5-0}
yields
$$
\E\,I_\zs{n}(f)=(\varrho^2_\zs{1}+\varrho^2_\zs{2}\lambda)||f||^{2}\,.
$$
In order to meet the condition \eqref{sec:In.5-1}
 the coefficients $\varrho_\zs{1}$, $\varrho_\zs{2}$
and the intensity  $\lambda>0$ must satisfy the inequality
\begin{equation}\label{sec:Ex.2}
\varrho^{2}_\zs{1}+\varrho^{2}_\zs{2} \lambda
\le
\sigma^{*}\,.
\end{equation}
 Note that the coefficients $\varrho_\zs{1}$,
$\varrho_\zs{2}$ and the intensity $\lambda$
in \eqref{sec:In.5-1}
 as well as $\sigma^{*}$
may depend on $n$, i.e.
$\varrho_\zs{i}=\varrho_\zs{i}(n)$ and $\lambda=\lambda(n)$.

As is stated in (\cite{KoPe3}, Theorem 2.2),
the conditions $\C_\zs{1})$ and $\C_\zs{2})$ hold for the process
\eqref{sec:Ex.1} with
$\sigma=\sigma(Q)=\varrho^2_\zs{1}+\varrho^2_\zs{2}\lambda$ defined in \eqref{sec:In.6},
$c^{*}_\zs{1}(n)=0$ and
$$
c^{*}_\zs{2}(n)\le 4\sigma(\sigma+\varrho^{2}_\zs{2}\E\,Y^{4}_\zs{1})\,.
$$
Let now $\cQ^{*}_\zs{n}$ be the family of distributions of the processes
\eqref{sec:Ex.1} with the coefficients
satisfying the conditions \eqref{sec:Ex.2}
and
\begin{equation}\label{sec:Ex.3}
\varrho^{2}_\zs{2}\,\le\, \sqrt{l_\zs{n}}\,,
\end{equation}
where
the sequence $l_\zs{n}$ is taken from the definition of the set $\cP^{*}_\zs{n}$.
Note that the distribution $Q_\zs{0}$ belongs to $\cQ^{*}_\zs{n}$. One can obtain this distribution
putting in \eqref{sec:Ex.1} $\varrho_\zs{1}=\sqrt{\sigma^{*}}$ and
$\varrho_\zs{2}=0$.
It will be noted that $\cQ^{*}_\zs{n}\subset\cP^{*}_\zs{n}$ if
$$
4\sigma^{*}(
\sigma^{*}+\sqrt{l}_\zs{n}\E\,Y^{4}_\zs{1}
)
\le l_\zs{n}\,.
$$

\section{Upper bound}\label{sec:Up}
\subsection{Known smoothness}

First we suppose that the parameters $k\ge 1$, $r>0$ and $\sigma^{*}$ in
\eqref{sec:In.5-1} are known. Let the family of admissible
weighted least squares estimates
$(\wh{S}_\zs{\gamma})_\zs{\gamma\in\Gamma}$
for the unknown function $S\in W^{k}_\zs{r}$ be given
\eqref{sec:Or.2}, \eqref{sec:Or.2-1}. Consider the pair
$$
\alpha_\zs{0}=(k,t_\zs{0})
$$
where
 $t_\zs{0}=[\ov{r}_\zs{n}/\ve]\ve$,
$\ov{r}_\zs{n}=r/\sigma^{*}_\zs{n}$
and $\ve$ satisfies the conditions
in \eqref{sec:Or.1-1}.
Denote the corresponding weight sequence in $\Gamma$ as
\begin{equation}\label{sec:Up.1}
\gamma_\zs{0}=\gamma_\zs{\alpha_\zs{0}}\,.
\end{equation}

Note that for sufficiently large $n$ the parameter
$\alpha_\zs{0}$ belongs to the set \eqref{sec:Or.2}. In this
section we obtain the upper bound for the empiric squared error of
the estimator \eqref{sec:In.12}.

\begin{theorem}\label{Th.sec:Up.1}
The estimator $\wh{S}_\zs{\gamma_\zs{0}}$ satisfies the following asymptotic upper
bound
\begin{equation}\label{Sec:Up.2}
\limsup_{n\to\infty}\,n^{2k/(2k+1)}\,
\frac{1}{R^{*}_\zs{k,n}}\,
\sup_{S\in W^{k}_\zs{r}}\,
\cR^{*}_\zs{n}\,(\wh{S}_\zs{\gamma_\zs{0}},S)\,
\le 1\,.
\end{equation}
\end{theorem}
\noindent {\bf Proof.}
First by substituting the model \eqref{sec:In.1} in the definition of
$\wh{\theta}_\zs{j,n}$ in \eqref{sec:In.1-1}
we obtain
$$
\wh{\theta}_\zs{j,n}=
\theta_\zs{j}
+
\frac{1}{\sqrt{n}}\xi_\zs{j,n}\,,
$$
where the random variables $\xi_\zs{j,n}$ are defined in
\eqref{sec:In.6}. Therefore, by the definition of
the estimators $\wh{S}_\zs{\gamma}$ in \eqref{sec:In.1-1}
we get
$$
||\wh{S}_\zs{\gamma_\zs{0}}-S||^2
=\sum_{j=1}^{n}\,(1-\gamma_\zs{0}(j))^2\,\theta^2_\zs{j}-2M_\zs{n}
 +\, \sum_{j=1}^{n}\,\gamma_\zs{0}^2(j)\,\xi^2_\zs{j,n}
$$
with
$$
M_\zs{n}\,=\,\frac{1}{\sqrt{n}}
\sum_{j=1}^{n}\,(1\,-\,\gamma_\zs{0}(j))\,\gamma_\zs{0}(j)\,\theta_\zs{j}\,\xi_\zs{j,n}
\,.
$$
It should be observed that
$$
\E_\zs{Q,S}\,M_\zs{n}=0
$$
for any $Q\in\cQ^{*}_\zs{n}$. Further the condition \eqref{sec:In.5-1} implies
also the inequality
$\E_\zs{Q}\,\xi^{2}_\zs{j,n}\le \sigma^{*}_\zs{n}$ for each distribution
$Q\in\cQ^{*}_\zs{n}$. Thus,
\begin{equation}\label{sec:Up.3}
\cR^{*}_\zs{n}(\wh{S}_\zs{\gamma_\zs{0}},S)
\,\le\,
\sum_{j=\iota_\zs{0}}^{n}\,(1-\gamma_\zs{0}(j))^2\,\theta^2_\zs{j}+
\frac{\sigma^{*}_\zs{n}}{n}\sum_{j=1}^{n}\,\gamma_\zs{0}^2(j)
\end{equation}
where $\iota_\zs{0}=j_\zs{0}(\alpha_\zs{0})$. Denote
$$
\upsilon_\zs{n}= n^{2k/(2k+1)} \sup_\zs{j\ge
\iota_\zs{0}}(1-\gamma_\zs{0}(j))^2/a_\zs{j}\,,
$$
where $a_\zs{j}$ is the sequence as defined in \eqref{sec:In.14}. Using this sequence
we estimate the first summand in the right hand of \eqref{sec:Up.3}
as
$$
n^{2k/(2k+1)}\,
\sum_{j=\iota_\zs{0}}^{n}\,
(1-\gamma_\zs{0}(j))^2\,\theta^2_\zs{j}
\le
\upsilon_\zs{n}\,
\sum_\zs{j\ge 1}\,a_\zs{j}\,\theta^{2}_\zs{j}
\,.
$$
From here and
\eqref{sec:In.14} we obtain that for each $S\in W^{k}_\zs{r}$
$$
\Upsilon_\zs{1,n}(S)=\,n^{2k/(2k+1)}\,
\sum_{j=\iota_\zs{0}}^{n}\,(1-\gamma_\zs{0}(j))^2\,\theta^2_\zs{j}
\,\le\,\upsilon_\zs{n}\,r\,.
$$
Further we note that
$$
\limsup_{n\to\infty}\,\left(\ov{r}_\zs{n}\right)^{2k/(2k+1)}\,
\upsilon_\zs{n}
\le\,
\frac{1}{\pi^{2k}\left(\tau_\zs{k}\right)^{2k/(2k+1)}}
\,,
$$
where the coefficient $\tau_\zs{k}$ is given \eqref{sec:Or.2-1}.
Therefore, for any $\eta>0$ and sufficiently large $n\ge 1$
\begin{equation}\label{Up.4}
\sup_\zs{S\in W^{k}_\zs{r}}\,
\Upsilon_\zs{1,n}(S)
\,\le\,
(1+\eta)\,(\sigma^{*}_\zs{n})^{2k/(2k+1)}\,
\Upsilon^{*}_\zs{1}
\end{equation}
where
$$
\Upsilon^{*}_\zs{1}=
\frac{r^{1/(2k+1)}}{\pi^{2k}(\tau_\zs{k})^{2k/(2k+1)}}\,.
$$
To examine the second summand in the right hand of \eqref{Sec:Up.2} we set
$$
\Upsilon_\zs{2,n}=\frac{1}{n^{1/(2k+1)}}\sum_{j=1}^{n}\,\gamma_\zs{0}^2(j)\,.
$$
Since by the condition \eqref{sec:In.12-1}
$$
\lim_\zs{n\to\infty} \,\frac{t_\zs{0}}{\ov{r}_\zs{n}}=1\,,
$$
one gets
$$
\lim_{n\to\infty}\,
\frac{1}{(\ov{r}_\zs{n})^{1/(2k+1)}}\,
\Upsilon_\zs{2,n}=
\Upsilon^{*}_\zs{2}
\quad\mbox{with}\quad
\Upsilon^{*}_\zs{2}=
\frac{2(\tau_\zs{k})^{1/(2k+1)}\,k^2}{(k+1)(2k+1)}\,.
$$
Note that  by the definition \eqref{sec:Mr.2}
$$
(\sigma^{*}_\zs{n})^{2k/(2k+1)}\,
\Upsilon^{*}_\zs{1,n}
+
\sigma^{*}_\zs{n}(\ov{r}_\zs{n})^{1/(2k+1)}
\Upsilon^{*}_\zs{2}
=
R^{*}_\zs{k,n}\,.
$$
Therefore, for any $\eta>0$ and sufficiently large $n\ge 1$
$$
n^{2k/(2k+1)}\,\sup_\zs{S\in W^{k}_\zs{r}}
\,\cR^{*}_\zs{n}(\wh{S}_\zs{\gamma_\zs{0}},S)
\le\,
(1+\eta)R^{*}_\zs{k,n}\,.
$$
Hence Theorem~\ref{Th.sec:Up.1}.
\endproof

\subsection{Unknown smoothness}

Combining Theorem~\ref{Th.sec:Up.1} and Theorem~\ref{Th.sec:Or.1}
yields Theorem~\ref{Th.sec:Mr.1}.

\endproof

\section{Lower bound}\label{sec:Lo}

First we obtain the lower bound for the risk \eqref{sec:In.5}
in the case of
 "white noise" model \eqref{sec:In.1}, when
 $\xi_\zs{t}=\sqrt{\sigma^{*}}w_\zs{t}$.
As before let $Q_\zs{0}$ denote the distribution of
$(\xi_\zs{t})_\zs{0\le t\le n}$ in $\cD[0,n]$.

\begin{theorem}\label{Th.sec:Lo.1}
The  risk \eqref{sec:In.5} corresponding to the the distribution
$Q_\zs{0}$ in the model \eqref{sec:In.1} has the following
lower bound
 \begin{equation}\label{Sec:Lo.0}
\liminf_{n\to\infty}\,
n^{2k/(2k+1)}\,
\inf_{\wt{S}_\zs{n}\in\Pi_\zs{n}}\,
\frac{1}{R^{*}_\zs{k,n}}\,
\sup_\zs{S\in W^{k}_\zs{r}}\, \cR_\zs{0}(\wt{S}_\zs{n},S)\, \ge\,
1\,,
\end{equation}
where $\cR_\zs{0}(\cdot,\cdot)=\cR_\zs{Q_\zs{0}}(\cdot,\cdot)$.
\end{theorem}
\proof
The proof of this result proceeds along the lines of Theorem 4.2 from
\cite{GaPe3}. Let $V$ be a function from $\C^{\infty}(\bbr)$ such that
$V(x)\ge 0$, $\int^{1}_\zs{-1}V(x)\d x=1$ and $V(x)=0$ for $|x|\ge 1$.
For each $0<\eta<1$ we introduce a smoother indicator of the interval
$[-1+\eta,1-\eta]$ by the formula
$$
I_\zs{\eta}(x)
=\eta^{-1}\int_{\bbr}\Chi_\zs{(|u|\le 1-\eta)}G\left(\frac{u-x}{\eta}\right)\,\d u\,.
$$
It will be noted that  $I_\zs{\eta}\in \cC^{\infty}(\bbr)$, $0\le I_\zs{\eta}\le 1$
and for any $m\ge 1$ and positive constant $c>0$
\begin{equation}\label{Sec:Lo.1}
\lim_\zs{\eta\to 0}\,
\sup_\zs{\{f\,:\,|f|_\zs{*}\le c\}}\,
\left|
\int_{\bbr}f(x)I^m_{\eta}(x)\,\d x\,-\,\int_{-1}^{1}f(x)\,\d x
\right|
=0
\end{equation}
where $|f|_\zs{*}=\sup_\zs{-1\le x\le 1}|f(x)|$.
Further, we need the trigonometric basis
 in $\cL_\zs{2}[-1,1]$, that is
\begin{equation}\label{Sec:Lo.2}
e_\zs{1}(x)=1/\sqrt{2}\,,\  e_\zs{j}(x)\,=\,Tr_j(\pi[j/2] x)\,,\ j\ge 2\,.
\end{equation}

Now we will construct
of a family of approximation functions for a given regression
function $S$ following \cite{GaPe3}.
For fixed $0<\ve<1$ one chooses
the bandwidth function as
\begin{equation}\label{Sec:Lo.3}
h=h_\zs{n}=(\upsilon^{*}_\zs{\ve})^{\frac{1}{2k+1}}\,N_\zs{n}\,n^{-\frac{1}{2k+1}}
\end{equation}
with
$$
\upsilon^{*}_\zs{\ve}=\frac{\sigma^{*}_\zs{n}k\pi^{2k}}{(1-\ve)r2^{2k+1}(k+1)(2k+1)}
\quad\mbox{and}\quad
N_\zs{n}=\ln^4 n
$$
and considers the partition
of the interval
$[0,1]$
with the points $\wt{x}_\zs{m}=2hm$, $1\le m\le M$, where
$$
M=[1/(2h)]-1\,.
$$
For each interval $[\wt{x}_\zs{m}-h,\wt{x}_\zs{m}+h]$ we
specify the smoothed indicator as $I_\zs{\eta}(v_\zs{m}(x))$, where
$v_\zs{m}(x)=(x-\wt{x}_\zs{m})/h$. The approximation function
for $S(t)$ is given by
\begin{equation}\label{sec:Lo.4}
S_\zs{z,n}(x)=\sum_{m=1}^M\sum_{j=1}^N\,z_\zs{m,j}\,D_\zs{m,j}(x)\,,
\end{equation}
where
$z=(z_\zs{m,j})_\zs{ 1\le m\le M\,, 1\le j\le N}$ is an array of real numbers;
$$
D_\zs{m,j}(x)=e_j(v_m(x))I_\zs{\eta}\,(v_m(x))
$$
are orthogonal functions on $[0,1]$.

Note that the set $W^{k}_\zs{r}$ is a subset of the ball
$$
\B_\zs{r}=\{f\in\cL_\zs{2}[0,1]\,:\,||f||^{2}\le r\}\,.
$$
Now for a given estimate $\wt{S}_\zs{n}$ we
construct its projection
in $\cL_\zs{2}[0,1]$
 into $\B_\zs{r}$
$$
\wt{F}_\zs{n}:=Pr_\zs{\B_\zs{r}}(\wt{S}_\zs{n})\,.
$$
In view of the convexity of the set $\B_\zs{r}$
one has
$$
||\wt{S}_\zs{n}-S||^2\ge||\wt{F}_\zs{n}-S||^2
$$
for each $S\in W^{k}_\zs{r}\subset \B_\zs{r}$.

From here one gets the following inequalities
for the the risk \eqref{sec:In.5}
$$
\sup_\zs{S\in W^{k}_\zs{r}}\cR_\zs{0}(\wt{S}_\zs{n},S)\ge
\sup_\zs{S\in W^{k}_\zs{r}}\cR_\zs{0}(\wt{F}_\zs{n},S)
\ge
\sup_\zs{\{z\in \bbr^{d}\,:\,S_\zs{z,n}\in W^{k}_\zs{r}\}}\cR_\zs{0}(\wt{F}_\zs{n},S)\,,
$$
where $d=M N$.

In order to continue this chain of estimates we need to introduce a special prior
distribution on $\bbr^{d}$. Let
$\kappa=(\kappa_\zs{m,j})_\zs{1\le m\le M\,, 1\le j\le N}$
be a random array with the elements
\begin{equation}\label{sec:Lo.5}
\kappa_\zs{m,j}\,=\,t_\zs{m,j}\,\kappa^{*}_\zs{m,j}\,,
\end{equation}
where $\kappa^{*}_\zs{m,j}$ are i.i.d. gaussian $\cN(0,1)$ random variables and
the coefficients
$$
t_\zs{m,j}=\frac{\sqrt{\sigma^{*}_\zs{n}y^*_\zs{j}}}{\sqrt{nh}}\,.
$$
We choose the sequence $(y^*_\zs{j})_\zs{1\le j\le N}$ in the same way as
in \cite{GaPe3}( see (8.11)) , i.e.
$$
y^*_j\,
=\,N^{k}_\zs{n}\,j^{-k}-1
\,.
$$
We denote the distribution of $\kappa$ by $\mu_\zs{\kappa}$.
 We will consider it as a prior
distribution of the random parametric regression $S_\zs{\kappa,n}$
which is obtained from \eqref{sec:Lo.4} by replacing $z$ with $\kappa$.

Besides we introduce
\begin{equation}\label{Sec:Lo.7}
\Xi_\zs{n}=
\left\{
z\in\bbr^{d}\,:\,
\max_\zs{1\le m\le M}\,\max_\zs{1\le j\le N}\,
 \frac{|z_\zs{m,j}|}{t_\zs{m,j}}\le \ln n
 \right\}
\,.
\end{equation}
By making use of the distribution $\mu_\zs{\kappa}$, one obtains
$$
\sup_\zs{S\in W^{k}_\zs{r}}\,\cR_\zs{0}(\wt{S}_\zs{n},S) \ge\,
\int_{\{z\in\bbr^d\,:\,S_\zs{z,n}\in W^k_\zs{r}\}\cap\Xi_\zs{n}}\,
\E_\zs{Q_\zs{0},S_\zs{z,n}}||\wt{F}_\zs{n}-S_\zs{z,n}||^2\,\mu_{\kappa}(\d
z)\,.
$$
Further we introduce
 the Bayes risk as
$$
\wt{\cR}(\wt{F}_\zs{n})=
\int_\zs{\bbr^d}\cR_\zs{0}(\wt{F}_\zs{n},S_\zs{z,n})
\mu_\zs{\kappa}(\d z)
$$
and noting that $||\wt{F}_\zs{n}||^2\le r$ we come to the
inequality
\begin{equation}\label{Sec:Lo.12}
\sup_\zs{S\in W^{k}_\zs{r}}\,\cR_\zs{0}(\wt{S}_\zs{n},S)\,
\ge\,
\wt{\cR}(\wt{F}_\zs{n})\,-\,\varpi_\zs{n}
\end{equation}
where
$$
\varpi_\zs{n}=\E (\Chi_\zs{\{S_\zs{\kappa,n}\notin
W^k_\zs{r}\}}\,+\, \Chi_\zs{\Xi^c_\zs{n}})
(r+||S_\zs{\kappa,n}||^2)\,.
$$
 By Proposition~\ref{Pr.sec:Fa.1} from Appendix~\ref{sec:Fa} one has, for any $p>0$,
$$
\lim_\zs{n\to\infty}\,n^{p}\,\varpi_\zs{n}\,=0\,.
$$
Now we consider the first term in the right-hand side of \eqref{Sec:Lo.12}.
To obtain a lower bound for this term we use the $\cL_2[0,1]$-orthonormal
function family $(G_\zs{m,j})_\zs{1\le m\le M,1\le j\le N}$ which is defined as
$$
G_\zs{m,j}(x)=
\frac{1}{\sqrt{h}}e_j\left(v_m(x)\right)\Chi_\zs{\left(|v_m(x)|\le 1\right)}\,.
$$
We denote by $\wt{g}_\zs{m,j}$ and $g_\zs{m,j}(z)$ the Fourier coefficients
for functions $\wt{F}_\zs{n}$ and $S_z$, respectively, i.e.
$$
\wt{g}_\zs{m,j}=\int^{1}_\zs{0}\,\wt{F}_\zs{n}(x)G_\zs{m,j}(x)\d x
\quad\mbox{and}\quad
g_\zs{m,j}(z)=\int^{1}_\zs{0}\,S_\zs{z,n}(x)\,G_\zs{m,j}(x)\d x\,.
$$
Now it is easy to see that
$$
||\wt{F}_\zs{n}-S_\zs{z,n}||^2 \ge \sum_{m=1}^M\sum_{j=1}^N\,
(\wt{g}_\zs{m,j}-g_\zs{m,j}(z))^2\,.
$$
Let us introduce the functionals $K_\zs{j}(\cdot): \cL_\zs{1}[-1,1]\to\bbr$ as
$$
K_\zs{j}(f)=
\int^1_{-1}\,e^2_\zs{j}(v)\,f(v)\,\d v\,.
$$
In view of \eqref{sec:Lo.4}
we obtain that
$$
\frac{\partial }{\partial z_\zs{m,j}}
g_\zs{m,j}(z)=
\,
\int^{1}_\zs{0}\,D_\zs{m,j}(x)\,G_\zs{m,j}(x)\,\d x
=\,\sqrt{h}\,K_\zs{j}(I_\zs{\eta})\,.
$$
Now
Proposition~\ref{Pr.sec:Tr.1} implies
\begin{align*}
\wt{\cR}(\wt{F}_\zs{n})\,&\ge
\sum_{m=1}^M\sum_{j=1}^N\,
\int_\zs{\bbr^d}\,
\E_\zs{S_\zs{z,n}}\,
(\wt{g}_\zs{m,j}-g_\zs{m,j}(z))^2
\mu_\zs{\kappa}(\d z)\\
&\ge\,h
\sum_{m=1}^M\sum_{j=1}^N\,\frac{\sigma^{*}K^2_j(I_\zs{\eta})}
{K_\zs{j}(I^2_\zs{\eta})\,nh\,+\,t^{-2}_\zs{m,j}\sigma^{*}}\,.
\end{align*}
Therefore, taking into account the definition of the coefficients
$(t_\zs{m,j})$ in \eqref{sec:Lo.5}
we get
$$
\wt{\cR}(\wt{F}_\zs{n})\,\ge\,
\frac{\sigma^{*}}{2nh}\,
\sum_{j=1}^N\,\tau_\zs{j}(\eta,y^*_\zs{j})
$$
with
$$
\tau_\zs{j}(\eta,y)=\frac{K^2_\zs{j}(I_\zs{\eta})y}{K_\zs{j}(I^2_\zs{\eta})y+1}\,.
$$
Moreover, the limit equality \eqref{Sec:Lo.1} implies directly
$$
\lim_\zs{\eta\to 0}\,
\sup_\zs{j\ge 1}\,\sup_\zs{y\ge 0}\,
\left|
\frac{(y+1)\tau_\zs{j}(\eta,y)}{y}
-1
\right|
=0\,.
$$
Therefore, we can write that for any $\nu>0$
$$
\wt{\cR}(\wt{F}_\zs{n})\,\ge\,
\frac{\sigma^{*}}{2nh(1+\nu)}\,
\sum_{j=1}^N\,
\frac{y^*_\zs{j}}{y^*_\zs{j}+1}\,.
$$
It is easy to check directly that
$$
\lim_\zs{n\to\infty}\,
\frac{\sigma^{*}_\zs{n}}{2nh R^{*}_\zs{k,n}}\,
\sum_{j=1}^N\,
\frac{y^*_\zs{j}}{y^*_\zs{j}+1}
=
(1-\ve)^{\frac{1}{2k+1}}\,,
$$
where the coefficient $R^{*}_\zs{k,n}$ is defined in
\eqref{sec:Mr.1}. Therefore,  \eqref{Sec:Lo.12} implies  for any
$0<\ve<1$
$$
\liminf_\zs{T\to\infty}\inf_\zs{\wh{S}_\zs{n}}\,n^{\frac{2k}{2k+1}}\,
\frac{1}{R^{*}_\zs{k,n}}
\,\sup_\zs{S\in W^{k}_\zs{r}}\,\cR_\zs{0}(\wh{S}_\zs{n},S)\,
\ge\,
(1-\ve)^{\frac{1}{2k+1}}\,.
$$
Taking here limit as $\ve\to 0$ implies Theorem~\ref{Th.sec:Lo.1}.
\endproof

\medskip

\renewcommand{\theequation}{A.\arabic{equation}}
\renewcommand{\thetheorem}{A.\arabic{theorem}}
\renewcommand{\thesubsection}{A.\arabic{subsection}}
\section{Appendix}\label{sec:A}
\setcounter{equation}{0} \setcounter{theorem}{0}

\subsection{Properties of the parametric family \eqref{sec:Lo.4}}\label{sec:Fa}

In this subsection we consider the sequence of the random functions $S_\zs{\kappa,n}$
defined in \eqref{sec:Lo.4} corresponding to the random array
$\kappa=(\kappa_\zs{m,j})_\zs{1\le m\le M,1\le j\le N}$
given in \eqref{sec:Lo.5}.

\begin{proposition}\label{Pr.sec:Fa.1}
For any $p>0$
$$
\lim_\zs{n\to\infty}\,n^{p}\,\lim_\zs{n\to\infty}\,
\E\,||S_\zs{\kappa,n}||^{2}\left( \Chi_\zs{\{S_\zs{\kappa,n}
\notin W^{k}_\zs{r}\}} + \Chi_\zs{\Xi^{c}_\zs{n}} \right)=0\,.
$$
\end{proposition}

\noindent This proposition follows directly from Proposition 6.4 in \cite{GaPe4}.

\subsection{Lower bound for parametric "white noise" models.}\label{sec:Tr}

In this subsection we prove some version of the van Trees inequality from \cite{GiLe}
for the following model
\begin{equation}\label{sec:Tr.1}
\d y_\zs{t}=S(t,z)\d t+
\sqrt{\sigma^{*}}\,
\d w_\zs{t}\,,\quad 0\le t\le n\,,
\end{equation}
where $z=(z_\zs{1},\ldots,z_\zs{d})'$ is vector of unknown
parameters, $w=(w_\zs{t})_\zs{0\le t\le T}$ is a Winier process.
 We assume that the function
$S(t,z)$ is a linear function with respect to the
parameter $z$, i.e.
\begin{equation}\label{sec:Tr.2}
S(t,z)=\sum^{d}_\zs{j=1}\,z_\zs{j}\,S_\zs{j}(t)\,.
\end{equation}
Moreover, we assume that the functions $(S_\zs{j})_\zs{1\le j\le d}$
are continuous.

Let $\Phi$ be a prior density in $\bbr^d$ having
the following form:
$$
\Phi(z)=\Phi(z_\zs{1},\ldots,z_\zs{d})=\prod_{j=1}^d\varphi_\zs{j}(z_\zs{j})\,,
$$
where $\varphi_\zs{j}$ is some continuously differentiable density
in $\bbr$. Moreover, let $g(z)$ be a continuously differentiable
$\bbr^d\to \bbr$ function such that for each $1\le j\le d$
\begin{equation}\label{sec:Tr.3}
\lim_\zs{|z_\zs{j}|\to\infty}\,
g(z)\,\varphi_\zs{j}(z_\zs{j})=0 \quad\mbox{and}\quad
\int_\zs{\bbr^d}\,|g^{\prime}_\zs{j}(z)|\,\Phi(z)\,\d z <\infty\,,
\end{equation}
where
$$
g^{\prime}_\zs{j}(z)=\frac{\partial g(z)}{\partial z_\zs{j}}\,.
$$
Let now $\cX_\zs{n}=\C[0,T]$ and $\cB(\cX_\zs{n})$ be $\sigma$ - field generated by
cylindric sets in $\cX_\zs{n}$.

 For any $\cB(\cX_\zs{n})\bigotimes \cB(\bbr^d)$-
measurable integrable function $\xi=\xi(x,\theta)$ we denote
$$
\wt{\E}\xi=\int_{\bbr^d}\, \int_\zs{\cX}\xi(y,z)\,
\mu_\zs{z}(\d y)\,\Phi(z) \d z \,,
$$
where $\mu_\zs{z}$ is distribution of the process \eqref{sec:Tr.1}
in $\cX_\zs{n}$. Let now $\nu=\mu_\zs{0}$ be the distribution of the process
$(\sigma^{*}w_\zs{t})_\zs{0\le t\le n}$ in $\cX$. It is clear (see, for example \cite{LiSh}) that
$\mu_\zs{z}<<\nu$ for any $z\in\bbr^{d}$.
Therefore, we can use the measure $\nu$ as a dominated measure, i.e.
for the observations \eqref{sec:Tr.1} in $\cX_\zs{n}$ we use
the following likelihood function
\begin{equation}\label{sec:Tr.4}
f(y,z)=
\frac{\d \mu_\zs{z}}{\d \nu}
=
\exp\left\{
\int^{n}_\zs{0}\,\frac{S(t,z)}{\sqrt{\sigma^{*}}}\,\d y_\zs{t}
-
\,
\int^{n}_\zs{0}\,\frac{S^2(t,z)}{2\sigma^{*}}\,\d t
\right\}\,.
\end{equation}

\vspace{5mm}

\begin{proposition}\label{Pr.sec:Tr.1}
  For any square integrable function
$\wh{g}_\zs{n}$
 measurable with respect to
 $\sigma\{y_\zs{t}\,,\,0\le t\le n\}$
 and for any $1\le j\le d$
  the following inequality holds
\begin{equation}\label{sec:Lo.6}
\tilde{\E}(\wh{g}_\zs{n}-g(z))^2\ge \frac{\sigma^{*}B^2_\zs{j}} {
\int_0^n\,S^{2}_\zs{j}(t)\,\d t
+\sigma^{*}I_\zs{j}}\,,
\end{equation}
 where
$$
B_\zs{j}=\int_\zs{\bbr^d}\,g^{\prime}_\zs{j}(z)\, \Phi(z)\,\d z
 \quad\mbox{and}\quad
I_\zs{j}=\int_\zs{\bbr}\,\frac{\dot{\varphi}^2_j(z)}{\varphi_\zs{j}(z)}\,\d
z\,.
$$
\end{proposition}
{\bf  Proof.} First of all note that the density \eqref{sec:Tr.3}
is bounded with respect to $\theta_\zs{j}\in\bbr$ for any $1\le
j\le d$, i.e. for any $y=(y_\zs{t})_\zs{0\le t\le n}\in \cX$
$$
\limsup_\zs{|z_\zs{j}|\to\infty}\,f(y, z)\,<\,\infty\,.
$$
Therefore, putting
$$
\Psi_j=\Psi_\zs{j}(y, z)=
 \frac{\partial}{\partial\theta_j}\,\ln(f(y, z)\Phi(z))
$$
and taking into account condition \eqref{sec:Tr.3} by integration
by parts one gets
\begin{align*}
\wt{\E}\left((\wh{g}_\zs{T}-g(z))\Psi_j\right)
&=\int_{\bbr^{N}\times\bbr^d}\,(\wh{g}_\zs{T}(y)-g(z))\frac{\partial}{\partial z_\zs{j}}
\left(f(y, z)\Phi(z)\right)\d z\,\d\nu( y)\\[2mm]
&=\int_{\bbr^{N}\times\bbr^d} \, g^{\prime}_\zs{j}(z)\, f(y,
z)\Phi(z)\,\d z \,\d \nu(y)=B_\zs{j}\,.
\end{align*}
Now by the Bounyakovskii-Cauchy-Schwarz inequality we obtain the
following lower bound for the quiadratic risk
$$
\wt{\E}(\wh{g}_\zs{T}-g(z))^2\ge
\frac{B^2_\zs{j}}{\wt{\E}\Psi_j^2}\,.
$$
 Note that from \eqref{sec:Tr.4} it is easy to deduce that under the distribution
 $\mu_\zs{z}$
\begin{align*}
\frac{\partial}{\partial z_\zs{j}}\,\ln f(y, z) &=
\,
\int^{n}_\zs{0}\,\frac{S_\zs{j}(t)}{\sqrt{\sigma^{*}}}\,\d y_\zs{t}
-
\,
\int^{n}_\zs{0}\,\frac{S(t,z)S_\zs{j}(t)}{\sigma^{*}}\,\d t\\
&=
\int^{n}_\zs{0}\,\frac{S_\zs{j}(t)}{\sqrt{\sigma^{*}}}\,\d w_\zs{t}
 \,.
\end{align*}
This implies directly
$$
\E_\zs{z}\, \frac{\partial}{\partial z_j}\,\ln f(y, z)\,=0
$$
and
$$
\E_\zs{z}\, \left(
\frac{\partial}{\partial z_j}\,\ln f(y, z)\, \right)^2
=\frac{1}{\sigma^{*}}\,\int^{n}_\zs{0}\,S^{2}_\zs{j}(t)\,\d t
\,.
$$
Therefore,
$$
\wt{\E}\Psi_j^2
=
\frac{1}{\sigma^{*}}\,\int^{n}_\zs{0}\,S^{2}_\zs{j}(t)\,\d t
+
I_\zs{j}\,.
$$
Hence Proposition~\ref{Pr.sec:Tr.1}.
\endproof


\medskip

\begin{thebibliography}{100}









\bibitem{Br}
Brua, J. (2009)
Asymptotically efficient estimators for nonparametric
heteroscedastic regression models.
{\sl Stat. Methodol.} 6(1), p. 47-60.


\bibitem{FoPe}
Fourdrinier, D. and Pergamenshchikov, S. (2007)
Improved selection model method for the regression with dependent noise.
 {\sl Annals of the Institute of Statistical Mathematics}, {\bf 59} (3),
p. 435-464.

\bibitem{GaPe0}
Galtchouk, L. and Pergamenshchikov, S. (2006)
Asymptotically efficient estimates for non parametric regression
models.{\sl Statistics and Probability Letters}, v. 76 , 8, p. 852-860

\bibitem{GaPe1}
Galtchouk, L. and Pergamenshchikov, S. (2004)
Nonparametric sequential estimation of the drift in diffusion processes.
{\sl Mathematical Methods of Statistics}, {\bf 13}, 1, 25-49.


\bibitem{GaPe2}
Galtchouk, L. and Pergamenshchikov, S. (2009)
 Sharp non-asymptotic oracle inequalities
for nonparametric heteroscedastic regression models.
{\em Journal of Nonparametric Statistics}, 2009, 21, 1, p. 1-16

\bibitem{GaPe3}
Galtchouk, L. and Pergamenshchikov, S. (2009) Adaptive
asymptotically efficient estimation in heteroscedastic
nonparametric regression. {\em Journal of Korean Statistical Society},

{\em http://ees.elsivier.com/jkss}

\bibitem{GaPe4}
Galtchouk, L. and Pergamenshchikov, S. (2009) Adaptive
asymptotically efficient estimation in heteroscedastic
nonparametric regression via model selection.

{\em http://hal.archives-ouvertes.fr/hal-00326910/fr/}






\bibitem{GiLe}
Gill, R.D. and Levit, B.Y.
Application of the van Trees inequality: a Bayesian Cram\'er-Rao bound.
{\em Bernoulli}, {\bf 1} 59-79 (1995)





\bibitem{KoPe2}
Konev, V.V. and Pergamenshchikov, S.M. (2008) General model
selection estimation of a periodic regression with a Gaussian
noise. - {\sl Annals of the Institute of Statistical Mathematics}, 2008,
Available online at\\
{\sl http://dx.doi.org/10.1007/s10463-008-0193-1}

\bibitem{KoPe3}
Konev, V.V. and Pergamenshchikov, S.M. (2009) Nonparametric
estimation in a  semimartingale regression model. Part 1. Oracle
Inequalities. {\sl Vestnik Tomskogo Universiteta, Mathematics and
Mechanics.}

\bibitem{LiSh}
Liptser, R. Sh. and Shiryaev, A.N. (1977)
{\sl Statistics of Random Processes. I. General theory.}
New York : Springer.



\bibitem{Nu}
Nussbaum, M. (1985) Spline smoothing in regression models and
asymptotic efficiency in $\L_2$. {\em Ann. Statist.}, {\bf 13}, p.
984-997.



\bibitem{Pi}
Pinsker, M.S. (1981) Optimal filtration of square integrable
signals in gaussian white noise. {\sl Problems of Transimission
information}, {\bf 17} 120--133
\end{thebibliography}
\end{document}